\documentclass[12pt]{article}
\usepackage{amsmath}
\usepackage{amssymb}
\usepackage{vatola}

\def\q{\quad}
\def\qq{\qquad}

\def\mod#1{\ (\text{\rm mod}\ #1)}
\def\t{\text}
\def\f{\frac}
\def\e{\equiv}
\def\b{\binom}
\def\sls#1#2{(\f{#1}{#2})}
 \def\ls#1#2{\big(\f{#1}{#2}\big)}
\def\Ls#1#2{\Big(\f{#1}{#2}\Big)}

\let \pro=\proclaim
\let \endpro=\endproclaim

\begin{document}
 \centerline {\bf
Elliptic curves  and the residue-counts of $x^2+bx+c/x$ modulo $p$}
\par\q\newline
\centerline{Zhi-Hong Sun}\newline
\centerline{School of Mathematics
and Statistics}
\centerline{Huaiyin Normal University}
\centerline{Huaian, Jiangsu 223300, P.R. China} \centerline{Email:
zhsun@hytc.edu.cn} \centerline{Homepage:
http://maths.hytc.edu.cn/szh1.htm}
 \abstract{\par For any prime $p>3$ and rational $p$-integers $b,c$ with $c(b^3-27c)\not\equiv 0\pmod p$ let
$V_p(x^2+bx+\frac cx)$ be the residue-counts of $x^2+bx+\frac cx$ modulo $p$ as $x$ runs over $1,2,\ldots,p-1$. In this paper, we reveal the connection between $V_p(x^2+bx+\frac cx)$ and the number of points on certain elliptic curve over the field $\Bbb F_p$.
 \par\q
\newline MSC(2020): Primary 11A15, Secondary 11A07, 11E25, 11L10, 11G20
 \newline Keywords: cubic residue; cubic congruence; elliptic curve; Jacobsthal sum; binary  quadratic form}
 \endabstract

\section*{1. Introduction}
\par\q\
For an odd prime $p$ let $\Bbb Z_p=\{0,1,\ldots,p-1\}$ and $\Bbb Z_p^*=\{1,2,\ldots,p-1\}$,
and let $R_p$ be the set of those rational numbers whose denominators are not divisible by $p$. For an odd prime $p$ and $a\in R_p$ let $\sls ap$ be the Legendre symbol.
 For an integral polynomial $f(x)$ let $V_p(f(x))$ be the residue-counts of $f(x)$ modulo $p$. Namely,
$V_p(f(x))$ is the number of $r\in\Bbb Z_p$ such that $f(x)\e r\mod p$ is solvable. For $c\in R_p$ let $V_p(f(x)+\f cx)$ be the residue-counts of $f(x)+\f cx$ modulo $p$. Namely,
$V_p(f(x)+\f cx)$ is the number of $r\in\Bbb Z_p^*$ such that $f(x)\e r\mod p$ is solvable.
 \par Let $p>3$ be a prime. It is well known that $V_p(x^2)=\f{p+1}2$. In 1908, von Sterneck[10] showed that for $a_1,a_2,a_3\in\Bbb Z$ with
$a_1^2-3a_2\not\e 0\mod p$,
$$V_p(x^3+a_1x^2+a_2x+a_3)=\f{2p+\sls p3}3.$$
See also [4, Theorem 4.3].  Suppose  $a,b\in R_p$ with $p\nmid b$. In [6], using the results for cubic and quartic congruences in [5] the author obtained a formula for $V_p(x^4+ax^2+bx)$. In particular, for $p\e 2\mod 3$ we have $V_p(x^4+bx)=[\f{5p+7}8]$, where $[\alpha]$ is the greatest integer not exceeding $\alpha$.
\par For $N\in\{14,15,20,24\}$ and $|q|<1$ let $a_N(n)$ be given by
$$\align &q\prod_{k=1}^{\infty}(1-q^k)(1-q^{2k})(1-q^{7k})(1-q^{14k})
=\sum_{n=1}^{\infty}a_{14}(n)q^n,
\\&q\prod_{k=1}^{\infty}(1-q^k)(1-q^{3k})(1-q^{5k})(1-q^{15k})
=\sum_{n=1}^{\infty}a_{15}(n)q^n,
\\&q\prod_{k=1}^{\infty}(1-q^{2k})^2(1-q^{10k})^2
=\sum_{n=1}^{\infty}a_{20}(n)q^n,
\\&q\prod_{k=1}^{\infty}(1-q^{2k})(1-q^{4k})(1-q^{6k})(1-q^{12k})
=\sum_{n=1}^{\infty}a_{24}(n)q^n.\endalign$$
It is well known (see [3]) that $a_N(n)$ is multiplicative and concerned with newforms of weight $2$ with level $N$.

\par Let $p>3$ be a prime. In [8], Sun and Ye obtained a formula for $V_p(x^3+\f cx)$ for $c\in \Bbb Z$ with $p\nmid c$. In particular, for $p\e 3\mod 4$,
$$V_p\Big(x^3+\f cx\Big)=\cases
\f{5p+5}8+\f 18\big(\frac{-c}{p}\big)a_{24}(p) &\t{if $8\mid p-7$ and $\ls{3c}p=1$,}
\\\f{5p-3}8+\f 18\big(\frac{-c}{p}\big)a_{24}(p) &\t{otherwise.}
\endcases$$
They also showed that for $p\e 3\mod 4$,
$$V_p(x^4-4x^2+4x)=\f 18\big(5p+1+4\delta(p)-2a_{20}(p)\big),$$
where
$$\delta(p)=\cases 0&\t{if
$p\e 7,23\mod {40}$,}
\\1&\t{if $p\e 3,27,31,39\mod {40}$,}
\\2&\t{if $p\e 11,19\mod {40}$.}\endcases$$
\par Let $p$ be a prime of the form $3k+1$. It is well known that $p=A^2+3B^2=\f 14(L^2+27M^2)$ with $A,B,L,M\in\Bbb Z$ and $A\e L\e 1\mod 3$. Jacobi
proved the following remarkable congruences (see [2]):
$$A\e \f 12\b{\f{p-1}2}{\f{p-1}6}\mod p\q\t{and}\q L\e -\b{\f{2(p-1)}3}{\f{p-1}3}\mod p.$$
\par Let $p>3$ be a prime and $a\in R_p$ with $a\not\e 0\mod p$. In [9] the author proved that
$$V_p\Big(x^2+\f{2a}x\Big)=\cases\f{2p-1}3\q\qq\qq\ \,\t{if $p\e 2\mod 3$,}
\\\f{2p-1+2A}3\q\qq\t{if $p=3k+1=A^2+3B^2(A,B\in\Bbb Z),\ 3\mid A-1$}
 \\\qq\qq\qq\qq\q\t{and $a^{\f{p-1}3}\e 1\mod p$,}
\\\f{2p-1-A\pm 3B}3\q\t{if $p=3k+1=A^2+3B^2(A,B\in\Bbb Z),\ 3\mid A-1$} \\\qq\qq\qq\qq\q\t{and $a^{\f{p-1}3}\e \f {-1\mp A/B}2\mod p$}.
\endcases$$
In Section 2, we give a new, natural and simple proof of the above result in [9].
\par
Let $p>3$ be a prime
 and
$b,c\in R_p$ with $bc(b^3-27c)\not\e 0\mod p$. Motivated by the above work, in Section 3 we evaluate $V_p(x^2+bx+\f cx)$ and reveal the connection between $V_p(x^2+bx+\f cx)$ and certain elliptic curve over $\Bbb F_p$. For  $m,n\in\Bbb Z$ let $\#
E_p(x^3+mx+n)$ be the number of points on the elliptic curve
$E_p:\ y^2=x^3+mx+n$ over the field $\Bbb F_p$ of $p$ elements.
It is easy to see (see for example [6, (2.5)]) that
$$\#
E_p(x^3+mx+n)=p+1+\sum_{x=0}^{p-1}\Ls{x^3+mx+n}p.\tag 1.1$$
 In Section 3 we show that for $t\in R_p$ with $t\not\e 0,\f 49\mod p$,
$$\#E_p(x^3+(6t-3)x+3t^2-6t+2)=\Big(2\Ls 3p+1\Big)p+1-\Ls 3p-3\Ls 3pV_p\Big(x^2+6x+\f{18t}x\Big).$$
As applications, we obtain
$$\align&V_p\Big(x^2+6x+\f 9x\Big)
\\&=\cases \f{2p-1}3&\t{if $p\e 2\mod 3$,}
\\\f{2p-1-L}3&\t{if $3\mid p-1$ and so $4p=L^2+27M^2\;(L,M\in\Bbb Z)$ with $3\mid L-1$}
\endcases\endalign$$
 and
$$V_p\Big(x^2+x-\f 1x\Big)=\f {2p-1+a_{14}(p)}3=\f{p-1}2+\f 18N(1,1,7,7;p)
=p-T(1,1,7,7;p-2),$$
where $N(a,b,c,d;n)$ is the number of representations of $n=ax^2+by^2+cz^2+dw^2$ $(x,y,z,w\in\Bbb Z)$, and $T(a,b,c,d;n)$ is the number of representations of $n=a\frac{x(x+1)}2+b\frac{y(y+1)}2+c\frac {z(z+1)}2+d\frac {w(w+1)}2$ $(x,y,z,w\in\{0,1,2,\ldots\})$.

\par Let $p>3$ be a prime, and $t\in R_p$ with $t\not\e 0,\f 49\mod p$. In Section 4 we obtain a new formula for $\#E_p(x^3-3x+2-9t)$ involving $V_p(x^4-6tx^2+12t^2x)$ and $V_p\big(x^2+6x+\f{18t}x\big)$.
As a consequence, we obtain a formula for $a_{15}(p)$. See Theorems 4.1 and 4.2.

\section*{2. Evaluation of $V_p\big(x^2+\f {2a}x\big)$}

\pro{Lemma 2.1} Let $p>3$ be a prime, $p\e 2\mod 3$, $m\in \Bbb Z$ and $m\not\e 0\mod p$.
Then
$$\sum_{x=1}^{p-1}\Ls{x(x^3+m)}p=\Ls mp\sum_{x=1}^{p-1}\Ls{x^3+m}p=-1.$$
\endpro
\par{\it Proof.} Since $p\e 2\mod 3$, for fixed $c\in\Bbb Z$ the congruence $x^3\e c\mod p$ has a unique solution. Thus,
$$\Ls mp\sum_{x=1}^{p-1}\Ls{x^3+m}p=\Ls mp\sum_{y=1}^{p-1}\Ls{y+m}p=\Ls mp\sum_{y=0}^{p-1}\Ls{y+m}p-1=-1.$$
Let $m'\in\Bbb Z$ be such that $mm'\e 1\mod p$. Then clearly
$$\align \sum_{x=1}^{p-1}\Ls {x(x^3+m)}p
&=\Ls mp\sum_{x=1}^{p-1}\Ls {x(m'x^3+1)}p
=\Ls mp\sum_{x=1}^{p-1}\Ls {m'+\f 1{x^3}}p
\\&=\Ls{m'}p\sum_{y=1}^{p-1}\Ls{y^3+m'}p=-1.
\endalign$$
This proves the lemma.
\vskip0.2cm
\pro{Lemma 2.2 ([2, Theorem 6.2.10])} Let $p$ be a prime of the form $3k+1$ and so $p=A^2+3B^2$ with $A,B\in\Bbb Z$ and $A\e 1\mod 3$. Suppose $m\in \Bbb Z$ and $m\not\e 0\mod p$.
Then
$$\sum_{x=1}^{p-1}\Ls{x(x^3+m)}p
=\cases -1-2A&\t{if $m^{\f{p-1}3}\e 1\mod p$,}
\\-1+A\pm 3B&\t{if $m^{\f{p-1}3}\e \f {-1\pm A/B}2\mod p$.}
\endcases$$\endpro

\pro{Theorem 2.1 ([9,Corollary 2.1])} Let $p>3$ be a prime and $a\in R_p$ with $a\not\e 0\mod p$.
\par $(\t{\rm i})$ If $p\e 2\mod 3$, then
$V_p(x^2+\f {2a}x)=\f{2p-1}3$.
\par $(\t{\rm ii})$ If $p\e 1\mod 3$
and so $p=A^2+3B^2$ with $A,B\in\Bbb Z$ and $A\e 1\mod 3$,
then
$$V_p\Big(x^2+\f{2a}x\Big)=\cases
\f{2p-1+2A}3&\t{if $a^{\f{p-1}3}\e 1\mod p$,}
\\\f{2p-1-A\pm 3B}3&\t{if $a^{\f{p-1}3}\e \f {-1\mp A/B}2\mod p$}
\endcases$$
and so
$a$ is a cubic residue of $p$ if and only if
 $V_p(x^2+\f {2a}x)=\f{2p-1+2A}3$.
\endpro
{\it Proof.} Suppose $x,y\in\Bbb Z_p^*$ and $x^2+\f{2a}x\e y^2+\f{2a}y\mod p$.
Then $x^2-y^2\e 2a(\f 1y-\f 1x)\mod p$ and so $(x+y)(x-y)\e 2a\f{x-y}{xy}\mod p$. If $x\not\e y\mod p$, we have
$x+y\e \f{2a}{xy}\mod p$. Hence,
$$y^2+xy\e \f{2a}x\mod p\q\t{and so}\q \Big(y+\f x2\Big)^2\e \f{x^3+8a}{4x}\mod p.$$
For fixed $x\in\Bbb Z_p^*$ with $\ls{x(x^3+8a)}p=0$, $\big(y+\f x2\big)^2\e \f{x^3+8a}{4x}\e 0\mod p$ if and only if $y\e -\f x2\mod p$, thus there exists a unique $y\in\Bbb Z_p^*$ such that
$x\not\e y\mod p$ and  $x^2+\f{2a}x\e y^2+\f{2a}y\mod p$.
For fixed $x\in\Bbb Z_p^*$ with $\ls{x(x^3+8a)}p=-1$,
$\big(y+\f x2\big)^2\e \f{x^3+8a}{4x}\mod p$ is unsolvable and so
$x^2+\f{2a}x\not\e y^2+\f{2a}y\mod p$ for $x\not\e y\mod p$.
For fixed $x\in\Bbb Z_p^*$ with $\ls{x(x^3+8a)}p=1$, there are exactly two $y\in\Bbb Z_p^*$ such that
$\big(y+\f x2\big)^2\e \f{x^3+8a}{4x}\mod p.$ Since $(x+\f x2)^2\e \f{x^3+8a}{4x}\mod p\iff x^3\e a\mod p$, we see that for fixed $x\in\Bbb Z_p^*$ with $\ls{x(x^3+8a)}p=1$ and $x^3\not\e a\mod p$, there are exactly two $y\in\Bbb Z_p^*$ such that $x\not\e y\mod p$ and $x^2+\f{2a}x\e y^2+\f{2a}y\mod p$. For $x\in\Bbb Z_p^*$ with $x^3\e a\mod p$, we have
$\ls{x(x^3+8a)}p=\ls{x^3(x^3+8a)}p=\ls{9a^2}p=1$ and $(x+\f x2)^2\e \f{x^3+8a}{4x}\mod p$, thus there exists a unique $y\in\Bbb Z_p^*$ such that
$x\not\e y\mod p$ and  $x^2+\f{2a}x\e y^2+\f{2a}y\mod p$.
\par For $a\in\Bbb Z_p^*$ set
$$\delta_p(a)=\cases 1&\t{if $p\e 2\mod 3$,}
\\0&\t {if $3\mid p-1$ and $a$ is not a cubic residue of $p$,}
\\3 &\t {if $3\mid p-1$ and $a$ is a cubic residue of $p$.}
\endcases$$
Then
$$\sum\Sb x\in\Bbb Z_p^*\\ p\mid x^3-a\endSb 1=\delta_p(a)=\sum\Sb x\in\Bbb Z_p^*\\ p\mid x^3+8a\endSb 1.$$
Summarizing the above deduces that
$$\aligned V_p\Big(x^2+\f{2a}x\Big)&=\f 12\sum\Sb x\in\Bbb Z_p^*\\\ls{x(x^3+8a)}p=0\endSb 1+\sum\Sb x\in\Bbb Z_p^*\\\ls{x(x^3+8a)}p=-1\endSb 1+\f 13\sum\Sb x\in\Bbb Z_p^*,p\nmid x^3-a\\\ls{x(x^3+8a)}p=1\endSb 1+\f 12\sum\Sb x\in\Bbb Z_p^*\\ p\mid x^3-a\endSb 1
\\&=\sum_{x\in\Bbb Z_p^*} 1-\f 23\sum\Sb x\in\Bbb Z_p^*,p\nmid x^3-a\\\ls{x(x^3+8a)}p=1\endSb 1-\f 12\sum\Sb x\in\Bbb Z_p^*\\ p\mid x^3-a\endSb 1-\f 12\sum\Sb x\in\Bbb Z_p^*\\ p\mid x^3+8a\endSb 1
\\&=p-1
-\f 23\sum\Sb x\in\Bbb Z_p^*\\\ls{x(x^3+8a)}p=1\endSb 1+\f 16\sum\Sb x\in\Bbb Z_p^*\\ p\mid x^3-a\endSb 1-\f 12\sum\Sb x\in\Bbb Z_p^*\\ p\mid x^3+8a\endSb 1
\\&=p-1-\f 13\delta_p(a)-\f 23\sum\Sb x\in\Bbb Z_p^*\\\ls{x(x^3+8a)}p=1\endSb 1.
\endaligned$$
On the other hand,
$$\sum\Sb x\in\Bbb Z_p^*\\\ls{x(x^3+8a)}p=1\endSb 1+\sum\Sb x\in\Bbb Z_p^*\\\ls{x(x^3+8a)}p=-1\endSb 1=p-1-\sum\Sb x\in\Bbb Z_p^*\\ p\mid x^3+8a\endSb 1=p-1-\delta_p(a)$$
and
$$\sum\Sb x\in\Bbb Z_p^*\\\ls{x(x^3+8a)}p=1\endSb 1-\sum\Sb x\in\Bbb Z_p^*\\\ls{x(x^3+8a)}p=-1\endSb 1=\sum_{x\in\Bbb Z_p^*}\Ls{x(x^3+8a)}p.$$
Thus,
$$2\sum\Sb x\in\Bbb Z_p^*\\\ls{x(x^3+8a)}p=1\endSb 1=p-1-\delta_p(a)
+\sum_{x\in\Bbb Z_p^*}\Ls{x(x^3+8a)}p$$
and so
$$\align V_p\Big(x^2+\f{2a}x\Big)&=p-1-\f 13\delta_p(a)-\f 13\Big(p-1-\delta_p(a)
+\sum_{x\in\Bbb Z_p^*}\Ls{x(x^3+8a)}p\Big)
\\&=\f 23(p-1)-\f 13\sum_{x=1}^{p-1}\Ls{x(x^3+8a)}p.\endalign$$
If $p\e 2\mod 3$, from Lemma 2.1 we have $\sum_{x=1}^{p-1}\ls{x(x^3+8a)}p=-1$, thus
$$V_p\Big(x^2+\f{2a}x\Big)=\f 23(p-1)-\f 13(-1)=\f{2p-1}3.$$
If $p\e 1\mod 3$ and so $p=A^2+3B^2$ with $A,B\in\Bbb Z$ and $A\e 1\mod 3$,
from the above and Lemma 2.2 we derive that
$$\align &V_p\Big(x^2+\f{2a}x\Big)\\&=\f 23(p-1)-\f 13\sum_{x=1}^{p-1}\Ls{x(x^3+8a)}p
\\&=\cases \f 23(p-1)-\f 13(-1-2A)=\f{2p-1+2A}3&\t{if $a^{\f{p-1}3}\e 1\mod p$,}
\\\f 23(p-1)-\f 13(-1+A\pm 3B)=\f{2p-1-A\mp 3B}3
&\t{if $a^{\f{p-1}3}\e \f {-1\pm A/B}2\mod p$.}
\endcases\endalign$$
Since $p$ is a prime we see that $A^2\not=B^2$ and so $A\not=\pm B$. Hence, if $a$ is a cubic non-residue modulo $p$, we have
 $V_p(x^2+\f {2a}x)=\f{2p-1+2A-3A\pm 3B}3\not=\f{2p-1+2A}3$.
 Thus $a$ is a cubic residue of $p$ if and only if
 $V_p(x^2+\f {2a}x)=\f{2p-1+2A}3$.
The proof is now complete.

\vskip 0.2cm

\section*{3. Evaluation of $V_p\big(x^2+bx+\f c{x}\big)$}

\par\q\  For a prime $p>3$ and integral polynomial $f(x)$ let $N_p(f(x))$ be the number of solutions to the congruence $f(x)\e 0\mod p$.
\pro{Theorem 3.1} Let $p>3$ be a prime, $b,c,t\in R_p$, $bc(b^3-27c)\not\e 0\mod p$ and $t\e \f{12c}{b^3}\mod p$. Then
$$3V_p\Big(x^2+bx+\f c{x}\Big)=2p-1-\Ls 3p\sum_{x=0}^{p-1}\Ls{x^3+(6t-3)x+3t^2-6t+2}p$$
and so for $t\not\e 0,\f 49\mod p$,
$$\#E_p(x^3+(6t-3)x+3t^2-6t+2)=\Big(2\Ls 3p+1\Big)p+1-\Ls 3p-3\Ls 3pV_p\Big(x^2+6x+\f{18t}x\Big).$$
\endpro
{\it Proof.} Suppose $x,y\in\Bbb Z_p^*$ and $x^2+bx+\f{c}x\e y^2+by+\f{c}y\mod p$.
Then $x^2-y^2+b(x-y)\e c\f {x-y}{xy}\mod p$. If $x\not\e y\mod p$, we have
$x+y+b\e \f c{xy}\mod p$ and so
$$ \Big(y+\f {x+b}2\Big)^2\e \f cx+\f{(x+b)^2}4=\f{x(x+b)^2+4c}{4x}\mod p.$$
For fixed $x\in\Bbb Z_p^*$ with $\ls{x(x+b)^2+4c)}p=0$, $\big(y+\f {x+b}2\big)^2\e \f{x(x+b)^2+4c}{4x}\e 0\mod p$ if and only if $y\e -\f {x+b}2\mod p$. Since $b^3\not\e 27c\mod p$ and $x(x+b)^2+4c)\e 0\mod p$, we have $x\not\e -\f {x+b}2\mod p$ and so there exists a unique $y\in\Bbb Z_p^*$ such that
$x\not\e y\mod p$ and  $x^2+bx+\f{c}x\e y^2+by+\f{c}y\mod p$.
For fixed $x\in\Bbb Z_p^*$ with $\ls{x(x(x+b)^2+4c)}p=-1$,
$\big(y+\f {x+b}2\big)^2\e \f{x(x+b)^2+4c}{4x}\mod p$ is unsolvable and so
$x^2+bx+\f{c}x\not\e y^2+by+\f{c}y\mod p$ for $x\not\e y\mod p$.
For fixed $x\in\Bbb Z_p^*$ with $\ls{x(x(x+b)^2+4c)}p=1$, there are exactly two $y\in\Bbb Z_p^*$ such that
$\big(y+\f {x+b}2\big)^2\e \f{x(x+b)^2+4c}{4x}\mod p.$ Since
$$\Big(x+\f {x+b}2\Big)^2\e \f{x(x+b)^2+4c}{4x}\mod p\iff 2x^3+ bx^2-c\e 0\mod p,$$
 for fixed $x\in\Bbb Z_p^*$ with $\ls{x(x(x+b)^2+4c)}p=1$ and $2x^3+bx^2\not\e c\mod p$, there are exactly two $y\in\Bbb Z_p^*$ such that $x\not\e y\mod p$ and $x^2+\f bx+\f{c}x\e y^2+by+\f cy\mod p$, and
 for fixed $x\in\Bbb Z_p^*$ with $\ls{x(x(x+b)^2+4c)}p=1$ and $2x^3+bx^2\e c\mod p$, there exists a unique $y\in\Bbb Z_p^*$ such that $x\not\e y\mod p$ and $x^2+\f bx+\f{c}x\e y^2+by+\f cy\mod p$.
 \par From the above we deduce that
 $$\aligned V_p\Big(x^2+bx+\f{c}x\Big)&=\f 12\sum\Sb x\in\Bbb Z_p^*\\x(x+b)^2+4c\e 0\mod p\endSb 1+\sum\Sb x\in\Bbb Z_p^*\\\sls{x(x(x+b)^2+4c)}p=-1\endSb 1\\&\q+\f 13\sum\Sb x\in\Bbb Z_p^*,\ p\nmid (2x^3+bx^2-c)\\\sls{x(x(x+b)^2+4c)}p=1\endSb 1+\f 12\sum\Sb x\in\Bbb Z_p^*,\ p\mid (2x^3+bx^2-c)\\\sls{x(x(x+b)^2+4c)}p=1\endSb 1
 \\&=\f 12\sum\Sb x\in\Bbb Z_p^*\\x(x+b)^2+4c\e 0\mod p\endSb 1+\sum\Sb x\in\Bbb Z_p^*\\\sls{x(x(x+b)^2+4c)}p=-1\endSb 1\\&\q+\f 13\sum\Sb x\in\Bbb Z_p^*\\\sls{x(x(x+b)^2+4c)}p=1\endSb 1+\f 16\sum\Sb x\in\Bbb Z_p^*,\ p\mid (2x^3+bx^2-c)\\\sls{x(x(x+b)^2+4c)}p=1\endSb 1
 \\&=p-1-\f 12N_p(x(x+b)^2+4c)
 -\f 23\sum\Sb x\in\Bbb Z_p^*\\\sls{x(x(x+b)^2+4c)}p=1\endSb 1
 \\&\q+\f 16\sum\Sb x\in\Bbb Z_p^*,\ p\mid (2x^3+bx^2-c)\\\sls{x(x(x+b)^2+4c)}p=1\endSb 1.
 \endaligned$$
Note that
$$\sum\Sb x\in\Bbb Z_p^*\\\sls{x(x(x+b)^2+4c)}p=1\endSb 1+\sum\Sb x\in\Bbb Z_p^*\\\sls{x(x(x+b)^2+4c)}p=-1\endSb 1=p-1-\sum\Sb x\in\Bbb Z_p^*\\ p\mid (x(x+b)^2+4c)\endSb 1=p-1-N_p(x(x+b)^2+4c)$$
and
$$\sum\Sb x\in\Bbb Z_p^*\\\sls{x(x(x+b)^2+4c)}p=1\endSb 1-\sum\Sb x\in\Bbb Z_p^*\\\sls{x(x(x+b)^2+4c)}p=-1\endSb 1=\sum_{x\in\Bbb Z_p^*}\Ls{x(x(x+b)^2+4c)}p.$$
We have
$$2\sum\Sb x\in\Bbb Z_p^*\\\sls{x(x(x+b)^2+4c)}p=1\endSb 1=p-1-N_p(x(x+b)^2+4c)
+\sum_{x\in\Bbb Z_p^*}\Ls{x(x(x+b)^2+4c)}p$$
and so
$$\aligned V_p\Big(x^2+bx+\f cx\Big)&=p-1-\f 12N_p(x(x+b)^2+4c)+\f 16\sum\Sb x\in\Bbb Z_p^*, p\mid (2x^3+bx^2-c)\\\sls{x(x(x+b)^2+4c)}p=1\endSb 1
\\&\q-\f 13\Big(p-1-N_p(x(x+b)^2+4c)+\sum_{x\in\Bbb Z_p^*}\Ls{x(x(x+b)^2+4c)}p\Big)
\\&=\f 23(p-1)-\f 16N_p(x(x+b)^2+4c)+\f 16\sum\Sb x\in\Bbb Z_p^*, p\mid (2x^3+bx^2-c)\\\sls{x(x(x+b)^2+4c)}p=1\endSb 1\\&\q-\f 13\sum_{x\in\Bbb Z_p^*}\Ls{x(x(x+b)^2+4c)}p.\endaligned$$
It is clear that
$$\align &\sum_{x\in\Bbb Z_p^*}\Ls{x(x(x+b)^2+4c)}p
\\&=\sum_{x\in\Bbb Z_p^*}\Ls{\f 1x(\f 1x(\f 1x+b)^2+4c)}p
=\sum_{x\in\Bbb Z_p^*}\Ls{\f 1{x^4}((1+bx)^2+4cx^3)}p
\\&=\sum_{x\in\Bbb Z_p^*}\Ls{4cx^3+(1+bx)^2}p=\sum_{x\in\Bbb Z_p^*}\Ls{4c^3x^3+(c+bcx)^2}p
\\&=\sum_{x\in\Bbb Z_p^*}\Ls{4x^3+(c+bx)^2}p
=\sum_{x=0}^{p-1}\Ls{4x^3+b^2x^2+2bcx+c^2}p-1
\\&=\sum_{x=0}^{p-1}\Ls{4(x-\f {b^2}{12})^3+b^2(x-\f {b^2}{12})^2+2bc(x-\f {b^2}{12})+c^2}p-1
\\&=\sum_{x=0}^{p-1}\Ls{4x^3+(2bc-\f {b^4}{12})x+c^2-\f{b^3c}6+\f{b^6}{216}}p-1
\\&=\Ls 2p\sum_{x=0}^{p-1}\Ls{(2x)^3+(2bc-\f {b^4}{12})\cdot 2x+2c^2-\f{b^3c}3+\f{b^6}{108}}p-1
\\&=\Ls 2p\sum_{x=0}^{p-1}\Ls{x^3+(2bc-\f {b^4}{12})x+2c^2-\f{b^3c}3+\f{b^6}{108}}p-1
\\&=\Ls 2p\sum_{x=0}^{p-1}\Ls{\sls x6^3+(2bc-\f {b^4}{12})\f x6+2c^2-\f{b^3c}3+\f{b^6}{108}}p-1
\\&=\Ls 3p\sum_{x=0}^{p-1}\Ls{x^3-3b(b^3-24c)x+2(b^6-36b^3c+216c^2)}p-1
.\endalign$$
If $2x^3+bx^2-c\e 0\mod p$, then clearly
$$\align x(x(x+b)^2+4c)&\e x(x(x+b)^2+8x^3+4bx^2)=x^2((x+b)^2+8x^2+4bx)
\\&=x^2(9x^2+6bx+b^2)=x^2(3x+b)^2\mod p.\endalign$$
Since $b^3\not\e 27c\mod p$, we see that $2(-\f b3)^3+b(-\f b3)^2-c\not\e 0\mod p$ and so $x\not\e -\f b3\mod p$. Also, $x\not\e 0\mod p$ since $c\not\e 0\mod p$.  Hence $\ls {x(x(x+b)^2+4c)}p=1$ and so
$$\aligned\sum\Sb x\in\Bbb Z_p^*,\ p\mid (2x^3+bx^2-c)\\\sls{x(x(x+b)^2+4c)}p=1\endSb 1
&=\sum_{x\in\Bbb Z_p^*,\ p\mid (2x^3+bx^2-c)} 1=N_p(2x^3+bx^2-c)
\\&=N_p((2x)^3+b(2x)^2-4c)=N_p(x^3+bx^2-4c)
\\&=N_p\Big(\big(x-\f b3\big)^3+b\big(x-\f b3\big)^2-4c\Big)
=N_p\Big(x^3-\f{b^2}3x+\f{2b^3}{27}-4c\Big)
\\&=N_p\Big(\Ls x3^3-\f{b^2}3\cdot \f x3+\f{2b^3}{27}-4c\Big)
=N_p(x^3-3b^2x+2b^3-108c).\endaligned$$
On the other hand,
$$\align &N_p(x(x+b)^2+4c)
\\&=N_p(x^3+2bx^2+b^2x+4c)=N_p\Big(\big(x-\f{2b}3\big)^3+2b\big(x-\f{2b}3\big)^2
+b^2\big(x-\f{2b}3\big)+4c\Big)
\\&=N_p\Big(x^3-\f{b^2}3x-\f 2{27}b^3+4c\Big)=N_p\Big(\big(-\f x3\big)^3-\f{b^2}3\big(-\f x3\big)-\f 2{27}b^3+4c\Big)
\\&=N_p(x^3-3b^2x+2b^3-108c).\endalign$$
\par Summarizing the above we derive that
$$\align &V_p\Big(x^2+bx+\f cx\Big)
\\& =\f 23(p-1)-\f 16N_p(x^3-3b^2x+2b^3-108c)+\f 16N_p(x^3-3b^2x+2b^3-108c)
\\&\q-\f 13\sum_{x\in\Bbb Z_p^*}\Ls{x(x(x+b)^2+4c)}p
\\&=\f 23(p-1)-\f 13\Big(\Ls 3p\sum_{x=0}^{p-1}\Ls{x^3-3b(b^3-24c)x+2(b^6-36b^3c+216c^2)}p-1\Big)
\\&=\f{2p-1}3-\f 13\Ls 3p\sum_{x=0}^{p-1}\Ls{x^3-3b(b^3-24c)x+2(b^6-36b^3c+216c^2)}p
.\endalign$$
Since $b\not\e 0\mod p$ and $12c\e b^3t\mod p$, we have
$$\align &V_p\Big(x^2+bx+\f cx\Big)
\\&=\f{2p-1}3-\f 13\Ls 3p\sum_{x=0}^{p-1}\Ls{x^3-3b(b^3-24c)x+2(b^6-36b^3c+216c^2)}p
\\&=\f{2p-1}3-\f 13\Ls 3p\sum_{x=0}^{p-1}\Ls{x^3-3b(b^3-2b^3t)x+2b^6-6b^6t+3b^6t^2)}p
\\&=\f{2p-1}3-\f 13\Ls 3p\sum_{x=0}^{p-1}\Ls{(b^2x)^3-3b^4(1-2t)b^2x+b^6(2-6t+3t^2)}p
\\&=\f{2p-1}3-\f 13\Ls 3p\sum_{x=0}^{p-1}\Ls{x^3+(6t-3)x+3t^2-6t+2}p
.\endalign$$
Hence, for $b=6$ and $c=18t$ we have
$$\Ls 3p\sum_{x=0}^{p-1}\Ls{x^3+(6t-3)x+3t^2-6t+2}p
=2p-1-3V_p\Big(x^2+6x+\f{18t}x\Big).$$
Now applying (1.1) deduces the remaining result.
\par\q

\pro{Theorem 3.2} Let $p>3$ be a prime. Then
$$\align&V_p\Big(x^2+6x+\f 9x\Big)
\\&=\cases \f{2p-1}3&\t{if $p\e 2\mod 3$,}
\\\f{2p-1-L}3&\t{if $3\mid p-1$ and so $4p=L^2+27M^2\;(L,M\in\Bbb Z)$ with $3\mid L-1$.}
\endcases\endalign$$
\endpro
{\it Proof.} Taking $b=6,\ c=9$ and $t=\f 12$ in Theorem 3.1 gives
$$\align 3V_p\Big(x^2+6x+\f 9x\Big)&=2p-1-\Ls 3p\sum_{x=0}^{p-1}\Ls{x^3-\f 14}p
\\&=2p-1-\Ls{-3}p-\Ls{-3}p\sum_{x=1}^{p-1}\Ls{1-4x^3}p
\\&=2p-1-\Ls{-3}p-\Ls{-3}p\sum_{x=1}^{p-1}\Ls{\f 1x(\f 1{x^3}-4)}p
\\&=2p-1-\Ls{-3}p-\Ls{-3}p\sum_{x=1}^{p-1}\Ls {x(x^3-4)}p
.\endalign$$
We first assume that $p\e 2\mod 3$. By Lemma 2.1, $\sum_{x=1}^{p-1}\ls
{x(x^3-4)}p=-1$. Thus,
$$V_p\Big(x^2+6x+\f 9x\Big)=\f 13\Big(2p-1-\Ls{-3}p+\Ls{-3}p\Big)=\f{2p-1}3.$$
\par Now assume that $p\e 1\mod 3$, $p=A^2+3B^2\; (A,B\in\Bbb Z)$, $4p=L^2+27M^2\;(L,M\in\Bbb Z)$ and $A\e L\e 1\mod 3$. In view of Lemma 2.2, we have
$$\align\sum_{x=1}^{p-1}\Ls {x(x^3-4)}p
&=\cases -1-2A&\t{if $(-4)^{\f{p-1}3}\e 1\mod p$,}
\\-1+A\pm 3B&\t{if $(-4)^{\f{p-1}3}\e\f{-1\pm A/B}2\mod p$}
\endcases
\\&=\cases -1-2A&\t{if $2^{\f{p-1}3}\e 1\mod p$,}
\\-1+A\pm 3B&\t{if $2^{\f{p-1}3}\e\f{-1\mp A/B}2\mod p$.}
\endcases
\endalign$$
From [6, (2.10)-(2.12)],
$$L=\cases -2A &\t{if $3\mid B$ and so $2^{\f{p-1}3}\e 1\mod p$,}
\\A+3B&\t{if $3\mid B-1$ and so $2^{\f{p-1}3}\e \f 12(-1-\f AB)\mod p$.}
\endcases$$
Hence,
$$\sum_{x=1}^{p-1}\Ls {x(x^3-4)}p
=\cases -1-2A&\t{if $3\mid B$,}
\\-1+A+3B&\t{if $3\mid B-1$}\endcases
=L-1.$$
Therefore,
$$V_p\Big(x^2+6x+\f 9x\Big)=
\f 13\Big(2p-1-1-\sum_{x=1}^{p-1}\Ls {x(x^3-4)}p\Big)
=\f{2p-1-L}3.$$
This completes the proof.
\par\q

\pro{Theorem 3.3} Let $p$ be a prime such that $p\not=2,3,7$. Then
$$\align &a_{14}(p)=3V_p\Big(x^2+x-\f 1x\Big)-2p+1,
\\&N(1,1,7,7;p)=8\Big(V_p\Big(x^2+x-\f 1x\Big)-\f{p-1}2\Big)
\\&T(1,1,7,7;p-2)=p-V_p\Big(x^2+x-\f 1x\Big).
\endalign$$
\endpro
{\it Proof.} Taking $b=6,\ c=-216$ and $t=-12$ in Theorem 3.1 gives
$$\align 3V_p\Big(x^2+6x-\f{216}x\Big)&=2p-1-\Ls 3p\sum_{x=0}^{p-1}
\Ls{x^3-75x+506}p
\\&=2p-1-\Ls 3p\sum_{x=0}^{p-1}
\Ls{(-x)^3-75(-x)+506}p
\\&=2p-1-\Ls{-3}p\Ls{x^3-75x-506}p.\endalign$$
Clearly,
$$V_p\Big(x^2+6x-\f{216}x\Big)=V_p\Big((6x)^2+6\cdot 6x-\f{216}{6x}\Big)
=V_p\Big(x^2+x-\f 1x\Big).$$
By [8, Lemma 2.1],
$$a_{14}(p)=-\Ls{-3}p\Ls{x^3-75x-506}p.$$
Thus
$$a_{14}(p)=3V_p\Big(x^2+6x-\f{216}x\Big)-(2p-1)=3V_p\Big(x^2+x-\f 1x\Big)-2p+1.$$
\par By [7, Lemma 2.10],
$$N(1,1,7,7;p)=\f 43(p+1)+\f 83a_{14}(p).$$
Thus,
$$\align N(1,1,7,7;p)&=\f 43(p+1)+\f 83\Big(3V_p\Big(x^2+x-\f 1x\Big)-2p+1\Big)
\\&=8\Big(V_p\Big(x^2+x-\f 1x\Big)-\f{p-1}2\Big).\endalign$$
Since $a_{14}(n)$ is multiplicative, $a_{14}(2)=-1$ and $a_{14}(4)=1$, from the above and [7, Theorem 4.16] we deduce that
$$\align T(1,1,7,7;p-2)&=\f{p+1}3+\f 16(a_{14}(2p)-a_{14}(4p))
\\&=\f{p+1}3+\f 16(a_{14}(2)-a_{14}(4))a_{14}(p)
=\f{p+1}3-\f 13a_{14}(p)
\\&=\f 13\Big(p+1+2p-1-3V_p\Big(x^2+x-\f 1x\Big)\Big)=p-V_p\Big(x^2+x-\f 1x\Big).
\endalign$$
This completes the proof.

\section*{4. A formula for $\#E_p(x^3-3x+2-9t)$}

\pro{Lemma 4.1 ([6, Corollary 2.2])} Let $p>3$ be a prime, $k\in R_p$ and $k\not\e 0\mod p$.
Then
$$\align &V_p(x^4+2kx^2+4k^2x)\\&=\f 18\Big\{5p+2+(-1)^{\f{p-1}2}+4\delta(k,p)
+\Ls p3 \big\{\# E_p(x^3-(18k+3)x-27k^2-18k-2)\\&\q-\#
E_p(x^3-3k^2x+k^3(27k+2))\big\} \Big\},
\endalign$$
 where
$\delta(k,p)$ is given by
$$\delta(k,p)=\Big|\Big\{x\bigm| x^3+4kx+8k^2\e 0 \mod p,\ \Ls xp=-1,
\ x\in\Bbb Z_p^*\} \Big\}\Big|.$$\endpro
\pro{Lemma 4.2 ([6, Remark 2.1])} Let $p>3$ be a prime, $k\in R_p$ and $k\not\e 0\mod p$. Then
$$\#E_p(x^3-3k^2x+k^3(27k+2))=\Big(1-\Ls kp\Big)(p+1)+\Ls kp\#E_p(x^3-3x+27k+2).$$

\pro{Theorem 4.1}  Let $p>3$ be a prime, $t\in R_p$ and $t\not\e 0,\f 49\mod p$. Then
$$\align &\Ls tp\big(\#E_p(x^3-3x+2-9t)-p-1\big)\\&=7p+1+(-1)^{\f{p-1}2}+4\varepsilon_p(t)
-8V_p(x^4-6tx^2+12t^2x)-3V_p\Big(x^2+6x+\f{18t}x\Big),
\endalign$$
where
$$\varepsilon_p(t)=\Big|\Big\{x\bigm| x^3-3tx+3t^2\e 0\mod p,
\ \Ls {6x}p=-1,\ x\in\Bbb Z_p^*\Big\}\Big|.$$
\endpro
{\it Proof.} Taking $k=-\f t3$ in Lemma 4.1 and then applying Lemma 4.2 gives
$$\align &8V_p\Big(x^4-\f{2t}3x^2+\f{4t^2}9x\Big)
\\&=5p+2+(-1)^{\f{p-1}2}+4\delta\Big(-\f t3,p\Big)
+\Ls p3 \Big\{\# E_p(x^3+(6t-3)x-3t^2+6t-2)\\&\q-\Big(1-\Ls {-3t}p\Big)(p+1)-\Ls {-3t}p\#E_p(x^3-3x+2-9t)\Big\}.
\endalign$$
Note that
$$V_p\Big(x^4-\f{2t}3x^2+\f{4t^2}9x\Big)
=V_p\Big(\ls x3^4-\f{2t}3\ls x3^2+\f{4t^2}9\cdot \f x3\Big)
=V_p(x^4-6tx^2+12t^2x)$$
and
$$\align\delta\Big(-\f t3,p\Big)&=\Big|\Big\{x\bigm| x^3-\f{4t}3x+\f 89t^2\e 0\mod p,
\ \Ls xp=-1,\ x\in\Bbb Z_p^*\Big\}\Big|
\\&=\Big|\Big\{x\bigm| \ls {2x}3^3-\f{4t}3\cdot\f{2x}3+\f 89t^2\e 0\mod p,
\ \Ls {6x}p=-1,\ x\in\Bbb Z_p^*\Big\}\Big|
\\&=\Big|\Big\{x\bigm| x^3-3tx+3t^2\e 0\mod p,
\ \Ls {6x}p=-1,\ x\in\Bbb Z_p^*\Big\}\Big|=\varepsilon_p(t).\endalign$$
By (1.1) and Theorem 3.1,
$$\align&\#E_p(x^3+(6t-3)x-3t^2+6t-2)
\\&=p+1+\sum_{x=0}^{p-1}\Ls{x^3+(6t-3)x-3t^2+6t-2}p
\\&=p+1+\sum_{x=0}^{p-1}\Ls{(-x)^3+(6t-3)(-x)-3t^2+6t-2}p
\\&=p+1+\Ls{-1}p\sum_{x=0}^{p-1}\Ls{x^3+(6t-3)x+3t^2-6t+2}p
\\&=p+1+\Ls{-3}p\Big(2p-1-3V_p\Big(x^2+6x+\f{18t}x\Big)\Big).
\endalign$$
Now, putting all the above together we deduce that
$$\align &8V_p(x^4-6tx^2+12t^2x)-5p-2-(-1)^{\f{p-1}2}-4\varepsilon_p(t)
\\&=\Ls{-3}p\Big(p+1+\Ls{-3}p\Big(2p-1-3V_p\Big(x^2+6x+\f{18t}x\Big)
\\&\q-\Big(1-\Ls{-3t}p\Big)(p+1)-\Ls{-3t}p\#E_p(x^3-3x+2-9t)\Big)\Big)
\\&=\Ls{t}p(p+1)+2p-1-3V_p\Big(x^2+6x+\f{18t}x\Big)-\Ls tp\#E_p(x^3-3x+2-9t),
\endalign$$
which yields the result.
\par\q

\pro{Theorem 4.2} Let $p>5$ be a prime. Then
$$\align a_{15}(p)&=\f 13(2N(1,3,5,15;p)-p-1)
 \\&=\Ls p3\Big(8V_p(x^4-6x^2+2x)+3V_p\Big(x^2+x+\f 3x\Big)-7p-1-(-1)^{\f{p-1}2}-4\varepsilon_p\Big),\endalign$$
where
$$\varepsilon_p=\Big|\Big\{x\bigm| x^3-3x+18\e 0\mod p,
\ \Ls {x}p=-1,\ x\in\Bbb Z_p^*\Big\}\Big|.$$
\endpro
{\it Proof}. Taking $t=36$ in Theorem 4.1 gives that
$$\align &\#E_p(x^3-3x-322)-p-1
\\&=7p+1+(-1)^{\f{p-1}2}+4\varepsilon_p(36)
-8V_p(x^4-6^3x^2+2\cdot 6^3x)-3V_p\Big(x^2+6x+\f{18\cdot 36}x\Big).
\endalign$$
It is clear that
$$V_p\Big(x^4-6^3x^2+2\cdot 6^3x\Big)=V_p((6x)^4-6^3(6x)^2+2\cdot 6^3\cdot 6x)
=V_p(x^4-6x^2+2x)$$
and
$$V_p\Big(x^2+6x+\f{18\cdot 36}x\Big)=V_p\Big((6x)^2+6\cdot 6x+\f{18\cdot 36}{6x}\Big)=V_p\Big(x^2+x+\f 3x\Big).$$
Also,
$$\align \varepsilon_p(36)&=\Big|\Big\{x\bigm| x^3-3\cdot 36x+3\cdot 36^2\e 0\mod p,
\ \Ls {6x}p=-1,\ x\in\Bbb Z_p^*\Big\}\Big|
\\&=\Big|\Big\{x\bigm| (6x)^3-3\cdot 36\cdot 6x+3\cdot 36^2\e 0\mod p,
\ \Ls {x}p=-1,\ x\in\Bbb Z_p^*\Big\}\Big|
\\&=\Big|\Big\{x\bigm| x^3-3x+18\e 0\mod p,
\ \Ls {x}p=-1,\ x\in\Bbb Z_p^*\Big\}\Big|=\varepsilon_p.
\endalign$$
From the above and (1.1) we get
$$\align &\sum_{x=0}^{p-1}\Ls{x^3-3x-322}p
\\&=\#E_p(x^3-3x-322)-p-1
\\&=7p+1+(-1)^{\f{p-1}2}+4\varepsilon_p
-8V_p(x^4-6x^2+2x)-3V_p\Big(x^2+x+\f 3x\Big).
\endalign$$
By [8, Lemma 2.1],
$$\sum_{x=0}^{p-1}\Ls{x^3-3x-322}p=-\Ls p3a_{15}(p).$$
By [1, Theorem 3.3], $a_{15}(p)=\f 13(2N(1,3,5,15;p)-p-1).$
Thus the result follows.
$$\q$$
\par{\bf Acknowledgements.} The author was supported by the
National Natural Science Foundation of China (grant No. 12271200).

\end{document}